\documentclass[a4paper,11pt]{article}
\usepackage[margin=3cm]{geometry}
\usepackage{amsfonts,amsmath,amssymb,amscd,graphicx}
\usepackage[colorlinks]{hyperref}
\usepackage{mathrsfs}

\newtheorem{theorem}{Theorem}[section]
\newtheorem{corollary}[theorem]{Corollary}
\newtheorem{lemma}[theorem]{Lemma}

\newtheorem{proposition}[theorem]{Proposition}

\newtheorem{example}[theorem]{Example}
\numberwithin{equation}{section}

\newenvironment{preuve}[1][]
{\vskip 2mm  \noindent\emph{\bf Proof#1. }}{$\Box$ \vskip 2mm}

\newcommand{\beq}{\begin{eqnarray}}
\newcommand{\eeq}{\end{eqnarray}}
\newcommand{\beqe}{\begin{eqnarray*}}
	\newcommand{\eeqe}{\end{eqnarray*}}

\newcommand{\rk}[1]{{\color{red}#1}}

\DeclareMathOperator{\inj}{inj}
\DeclareMathOperator{\syst}{syst}

\DeclareMathOperator{\spec}{spec}
\DeclareMathOperator{\length}{length}
\DeclareMathOperator{\grass}{Grass}
\DeclareMathOperator{\dist}{dist}

\DeclareMathOperator{\diam}{diam}
\DeclareMathOperator{\ricci}{Ric}

\newcommand{\End}{\ensuremath \mathrm{End}}
\newcommand{\id}{\ensuremath \mathrm{Id}}

\newcommand{\Nn}{\mathbb{N}}
\newcommand{\R}{\mathbb{R}}

\newcommand{\C}{\mathbb{C}}

\newcommand{\ep}{\epsilon}

\let\epsilon=\varepsilon

\begin{document}

\title{Metric and spectral aspects of \\ random complex divisors}
\author{\sc Michele Ancona and Damien Gayet}
\maketitle
\begin{abstract}
For any integer $n\geq 2$, we prove that for any large enough integer $d$, with large probability the injectivity radius of a random degree $d$ complex hypersurface in $\C P^n$ 
is larger than 
$d^{-\frac{1}2(3n+2)}$. Here the hypersurface 
is endowed with the restriction of the ambient Fubini-Study metric, and the probability measure is
induced by the Fubini-Study $L^2$-Hermitian product on the space of homogeneous complex polynomials
of degree $d$ in $(n+1)$-variables. We also prove that with high probability,
the sectional curvatures of the random hypersurface are bounded by 
$d^{\frac{3}2(n+2)}$,
 and that its spectral gap is bounded below by 
$\exp(-d^{\frac{1}4(3n+15)})$. 
These results extend to random  submanifolds of higher codimension in any complex projective manifold. Independently, we prove that the diameter of a degree $d$ divisor is bounded by $Cd^3$, which generalizes and amends the bound given in~\cite{feng1999diameter} for planar curves.
\end{abstract}

\tableofcontents

\section{Introduction}
Smooth complex projective hypersurfaces in $\C P^n$ have the remarkable property that their
diffeomorphism type depends only on their degree. Moreover,  for $n\geq 2$ they
are all connected and for $n\geq 3$, they are all simply connected. For $n=2$, smooth complex curves of degree $d$ in $\C P^2$ are compact Riemann surfaces of 
genus $\frac{1}2(d-1)(d-2)$. When equipped with the
restriction of the Fubini-Study metric, these complex hypersurfaces become Riemannian manifolds, whose geometry strongly depends on the chosen hypersurfaces. In this paper, we study some of their metric properties when the hypersurface is taken at random. 
\subsection{Complex planar curves}
We begin to state the setting and the results in the special case where $n=2$, that is we consider the case of random complex curves in $\C P^2$. 
Let
\begin{equation}\label{mesure2} P= \sum_{i+j+k =d} a_{i,j,k}\sqrt{ \frac{(d+2)!}{2 i!j!k!}}
X_0^{i}X_1^jX_2^{k}, \end{equation}
be a degree $d$ random polynomial,
 where $(a_I)_{|I|=d}\in \C P^{N_d-1}$
are random coefficients chosen uniformly on the projective space $\C P^{N_d-1}$ equipped with its Fubini-Study metric (the quotient metric induced by the standard metric on $\C^{N_d}$) and $N_d= \dim \C^{hom}_d[X_0,X_1, X_2]$. We can also choose $(a_I)_I$ in the standard sphere $\mathbb S^{2N_d-1}\subset \C^{N_d}$. Remark that the Hermitian product associated with the orthonormal basis $\{\sqrt{ \frac{(d+2)!}{2 i!j!k!}}
X_0^{i}X_1^jX_2^{k}\}$ is very a natural Hermitian product on  the space of degree $d$ homogeneous polynomials: it is the only $U(3)$-invariant product on $\C^{hom}_d[X_0,X_1, X_2]$ and it also appears as the $L^2$-product associated with the Fubini-Study metric (see Example \ref{example}).

 The vanishing locus $Z(P)\subset \C P^2$ of $P$ is almost surely a smooth connected Riemann surface of genus $\frac{1}2(d-1)(d-2)$, which is naturally equipped with a Riemannian metric: the restriction $g_{FS|Z(P)}$ of the Fubini-Study metric $g_{FS}$ of $\C P^2$.
 As, almost surely, the topology of $Z(P)$ only depends on the degree $d$,  we can view $Z(P)$ as a \emph{fixed} genus $\frac{1}2(d-1)(d-2)$ surface equipped with a \emph{random} metric. Our first result concerns the systole, that is the smallest length of non-contractible loops.             
            \begin{theorem}\label{S2}
            	           	Under the setting given above,
            	$$  \ \mathbb \mu_d\left[ \syst(Z(P),g_{FS|Z(P)}) \geq 
 d^{-4}
            	\right]\underset{d\to \infty}{\to} 1.$$
            	Here, $\mu_d$ denotes the measure defined by~(\ref{mesure2}).
            \end{theorem}
          Thus, Theorem \ref{S2} gives a probabilistic lower bound for the systole of a complex plane curve. Remark that there are no non-trivial deterministic lower bounds for this, that is   for any  $d\geq 2$ one has $$ \inf_{P\in \C^{hom}_d[X_0,X_1,X_2]}       
                  \syst (Z(P),g_{FS|Z(P)}) =0.$$
                  Indeed, let $Q$ be the product of $d$ lines
            passing through the same point $M\in \C P^2$. Then,   any non-contractible curve in the smoothing of  $Z(Q)$ can be isotopied near $M$ to a curve with small length. More generally, if one consider a family of regular polynomials $\{P_t\}_{t\in [0,1[}$  that converges as $t\rightarrow 1$ to a  general point of the discriminant locus $\Delta_d\subset \C^{hom}_d[X_0,X_1,X_2]$ (that is, the locus of polynomials defining singular complex curves), then of can observe a non--contractible loop $\gamma_t$ of $Z(P_t)$ that shrinks more and more as $t\rightarrow 1$, until it becomes a double singular point at the limit $t=1$. We can then see that the length of the loop  $\gamma_t$ goes to zero as $t\rightarrow 1$.
 Remark that a theorem by Gromov~\cite[2.C]{gromov1992systoles} about the systole of compact real surfaces implies the following deterministic upper bound for the systole of a plane curve: 
           there exists $C>0$ such that for any $d$ and 
           any $\forall P\in \C^{hom}_d[X_0,X_1,X_2],$
           $$  
           \syst (Z(P),g_{FS|Z(P)})\leq C d^{-\frac12}\log d.$$            
      Moreover,  in \cite{gayet2022systoles}, it was proven that there exists $c>0$ such that for any $d$ large enough,
$$
\mu_d \left[\syst(Z(P),g_{FS|Z(P)})\leq d^{-\frac12}\right]\geq c.$$

           These probabilistic results for complex planar curves are inspired by similar ones for another model of random metrics over large genus compact surfaces. More precisely,
          in~\cite{mirzakhani13}, M. Mirzakhani studied probabilistic
           aspects of metric parameters of $(X,h)$, when 
           $(X,h)$ is taken at random in $\mathcal M_g$, the moduli
           space of hyperbolic genus $g$ compact Riemann surfaces. 
           This moduli space is equipped with a natural 
           symplectic form, the Weil-Petersson form, hence a volume form, for
           which $\mathcal M_g$ has a finite volume, and
           provides a natural probability measure $\mu_{WP,g}$ on it, see~\cite{mirzakhani13}.
           M. Mirzakhani proved  in this context~\cite[Theorem 4.2]{mirzakhani13}\label{mimi}:
           	there exists $\ep_0>0$ and $0<c<C$ such that 
           	for any $\ep\leq \ep_0 $ and every $g\geq 2$,
           	$$  \ c \ep^2 \leq \mu_{WP,g} \big[  \syst(X,h) \leq \ep
           	\big]\leq C \ep^2. $$ 
           
The Gaussian curvature $K$ of hyperbolic surfaces is constant and equal to $-1$. In our complex setting, the curvature is never constant, see~\cite{ness1977curvature}. Moreover, $K\leq 2$ and there always exists points on $Z(P)$ where $K=2$. Note also that by Gauss-Bonnet formula, 
$$ \int_{Z(P)}K(x)dx = 2-2g.$$ Since by Wirtinger theorem, the area of $Z(P)$ is $d$, the average over $Z(P)$ of the curvature is $-d$. In order to compare our results with the ones for hyperbolic surfaces, we should rescale $g_{FS}$ by $\sqrt d g_{FS}$ (the volume of $Z(P)$ for this metric is $d^2$). 
               \begin{theorem}\label{K2}
           	Under the setting given above,
      $$  \ \mathbb \mu_d\left[
           	\sup_{           			x\in Z(s)}|K(x)|
           	\leq d^{5} \right]\underset{d\to \infty}{\to}1.$$
           \end{theorem}
Notice that by~\cite{ness1977curvature}, for $d\geq 2$, 
           $$ \sup_{P\in \C^{hom}_d[X_0,X_1,X_2]} \sup_{x\in Z(P)}|K(x)|=+\infty.$$

           Finally, if $\lambda_1(Z(P),g_{FS|Z(P)})$ denotes the first positive eigenvalue of the Laplacian for $g_{FS|Z(P)},$ we obtain the following theorem:
 \begin{theorem}\label{L2}
           	Under the setting given above,
$$  \mathbb \mu_d\left[
	\lambda_1(Z(P),g_{FS|Z(P)})
	\geq \exp(-d^{\frac{11}2})\right]\underset{d\to \infty}{\to}1.$$
\end{theorem}
By arguments provided by Cheeger~\cite{cheeger}, one can prove that
$$ \inf_{P\in \C^{hom}_d[X_0,X_1,X_2]} 	\lambda_1(Z(P),g_{FS|Z(P)})= 0.$$
Moreover, by~\cite{bourguignon}, $$
	\sup_{P\in \C^{hom}_d[X_0,X_1,X_2]}	\lambda_1(Z(P),g_{FS|Z(P)})\leq 6.$$
	In the hyperbolic setting,  by~\cite{huber}, $$\limsup_{g\to \infty}\sup_{h\in \mathcal M_g}\lambda_1(X,h)\leq \frac14.$$
In~\cite{anantharaman2023friedmanramanujan}, it is proven that
in the hyperbolic setting, for any $\epsilon>0$,
$$ \mu_{WP,g}\left[\lambda_1(X,h)\geq \frac{2}9 -\epsilon\right] \underset{g\to \infty}{\to }1.$$
Table~\ref{tab1} (resp. Table~\ref{tab2}) sums up the differences and similarities between the deterministic (resp. probabilistic) 
results in the hyperbolic and complex projective settings. 
           \begin{figure}
           	$$\begin{array}{|l|l|l|}
           	\hline
           	&&\\
           	\text{Parameters} & \text{Hyperbolic surfaces} & \text{Planar algebraic curves}\\
           	\text{of the surfaces of genus }g & &g\sim \frac{1}2d^2 \\&& \\ \hline  
           	&&\\
           	\text{Dimension of the moduli space} & \underset{g\to\infty}{\sim}  3g & \underset{d\to\infty}{\sim}  g \\
           	           	           	&&\\
           	\text{Curvature}& -1 &\in ]-\infty, 2] \ \cite{ness1977curvature}\\
           	           	           	&&\\
           	\text{Volume} & \underset{g\to\infty}{\sim}  4\pi g &  \underset{g\to\infty}{\sim}  4\pi g\\
           	           	           	&&\\
           	\text{Diameter} & \in ]0,+\infty[ & \in [cd^{1/2}, C d^{5/2}]\\
           		&&(\text{Theorem }\ref{theorem: diameter})\\
           	&&\\

                      	\lambda_1 & \in ]0,\frac{1}4+o_g(1)[~\cite{huber} & \in ]0,\frac{6}d]~\cite{bourguignon}\\
           &&\\	
           	\hline
           	\end{array}$$
           	\caption{Deterministic parameters of the two different models of real surfaces, the Weil-Petersson one with hyperbolic surfaces, and the Fubini-Study model with complex algebraic curves equipped with the induced rescaled induced metric $\sqrt {2\pi d} g_{FS}$ on $\C P^2$.}\label{tab1}
           \end{figure}
           \begin{figure}
           	$$\begin{array}{|l|l|l|}
           	\hline
           	&&\\
           	\text{Parameters} & \text{Hyperbolic surfaces} & \text{Planar algebraic curves} \\
           	\text{of the surfaces}   & \text{measure } \mu_{WP,g} & \text{measure } \mu_d\\ 
           	\text{of genus }g &  &g\sim \frac{1}2d^2 \\ 
           	&&\\ \hline  
           	&&\\
           	\text{Curvature}
           	& -1  &\mu_d (\|K\|_\infty<d^{10})\underset{d\to\infty}{\to} 1  \\
           	&&(\text{Theorem }\ref{K2})\\
           	           	&&\\
           	\text{Diameter} & \mu_{WP, g} (\diam \geq 40\log g) \underset{g\to\infty}{\to}  0 \ \cite{mirzakhani13}& ? \\
           	&&\\
           	\text{Systole} & \mu_{WP, g} (\syst\leq \ep) \asymp \ep^2 \ \cite{mirzakhani13} & \mu_d (\syst\leq \ep) \geq \exp(- \frac{c}{\ep^{6}})~\cite{gayet2022systoles}\\
           	&& \mu_d (\syst \geq d^{-4})\underset{d\to\infty}{\to} 1\\
           	           	&&(\text{Theorem }\ref{S2})\\
           	           	&&\\
           	\lambda_1 & \mu_{WP, g} (\lambda_1\geq \frac29-\epsilon)\underset{d\to\infty}{\to} 1~\cite{anantharaman2023friedmanramanujan} &\mu_d(
           	\lambda_1 
           	\geq \exp(-d^{\frac{9}2}))\underset{d\to\infty}{\to} 1 \\
           		&&(\text{Theorem }\ref{L2})\\
           	&& \\
           	\hline
           	\end{array}$$
           	\caption{Statistics of some metric parameters. The complex algebraic curves are equipped with the induced rescaled induced metric $\sqrt {2\pi d} g_{FS}$ on $\C P^2$.}\label{tab2}
           \end{figure}    
  \subsection{The general K\"ahler setting}          
  The result of the previous section can be stated in a much more general setting that we now introduce.            
Let $X$ be a complex projective manifold of dimension $n\geq 1$ equipped with   a Hermitian ample holomorphic line bundle  $(L,h)\to X$  with positive curvature $\omega$, 
that is locally 
$$ \omega = \frac{1}{2i\pi}\partial \bar{\partial}\log \|s\|^2_h>0,$$ 
where $s$ is any local non vanishing holomorphic section of $L$. Let $g_\omega= \omega (\cdot, i\cdot)$ be the associated K\"ahler metric. Let $(E,h_E)\to X$ be a holomorphic vector bundle of rank equal to $r\geq 1$ equipped with a Hermitian metric $h_E$. 
The space $H^0(X,E\otimes L^{\otimes d})$ of holomorphic sections of $E\otimes L^{\otimes d}$ is non trivial for $d$ large enough, 
and can be equipped with the $L^2$ Hermitian product 
\begin{equation}\label{prod}
(s,t)\in (H^0(X,E\otimes L^d))^2\mapsto \langle s,t\rangle = \int_X \langle s(x),t(x)\rangle_{h_d}\frac{\omega^n}{n!},
\end{equation}
where 
$h_d:= h_E\otimes h^d.$
This product induces a Gaussian measure $\mu_d$ over $H^0(X,E\otimes L^d)$, 
that is for any Borelian $U\subset H^0(X,E\otimes L^d),$
\begin{equation}\label{mesure}
 \mu_d (U)=\int_{s\in U} e^{-\frac12 \|s\|^2} \frac{ds}{(2\pi)^{N_d}},
\end{equation}
where $N_d $ denotes the complex dimension of $H^0(X,E\otimes L^d)$ and $ds$ denotes the Lebesgue measure associated to the Hermitian product~(\ref{prod}).
If $(S_i)_{i\in \{1, \cdots, N_d\}}$
denotes an orthonormal basis of this space, then 
$$ s =\sum_{i=1}^{N_d} a_i S_i$$
follows the law $\mu_d$ if the random complexes $\sqrt 2 a_i$ are i.i.d standard complex Gaussians, 
that is $\Re a_i$ and $\Im a_i$ are independent centered Gaussian variables with variance equal to $1/2$.
Note that for any event depending only on the vanishing locus $Z(s)$ of $s\in H^0(X,E\otimes L^d)$, the probability measure $\mu_d$ can be replaced by the invariant measure over the unit sphere $$S^d = \mathbb S H^0(X,E\otimes L^d)$$ for the product~(\ref{prod}), or equivalently the Fubini-Study measure on the linear system $\mathbb P H^0(X,E\otimes L^d)$

\begin{example}\label{example} When $X= \C P^n$, $(L, h)=(\mathcal O(1),h_{FS})$ is the degree
1 holomorphic line bundle equipped with the standard Fubini-Study metric and $(E,h_E) = (\C P^n\times \C^r,h_0)$ is the trivial rank $r$ bundle equipped with the standard Hermitian product, then the vector
space $H^0(X,E\otimes L^d)$ is isomorphic to the space 
$ (\C_d^{hom}[Z_0 , \cdots , Z_n ])^r$ of $r$-uples of  degree $d$ homogeneous
 polynomials in $n + 1$ variables. In this case, if $(e_i)_{i\in \{1, \cdots, r\}}$ denotes the
 standard basis of $\C^r$, the family 
$$\left(\sqrt{ \frac{(n+d)!}{n!i_0!\cdots i_n!})}
Z_0^{i_0}\cdots Z_n^{i_n}\otimes e_i\right)_{i_0+\cdots +i_n =d, \ i\in \{1, \cdots, r\}}$$ 
is an orthonormal basis for the Hermitian metric~(\ref{prod}).
\end{example}

Let $\Delta_d\subset H^0(X,E\otimes L^d)$ be the discriminant subset, that is the set of sections $s$ such that there exists $x$ in $Z(s)$ where $\nabla s(x)$ is not onto. Recall that $\Delta_d$ is a complex subvariety of positive codimension, and that for any $s\in H^0(X,E\otimes L^d)\setminus \Delta_d$, the zero set $Z(s)\subset X$ is a compact smooth complex submanifold of $X$ of codimension $r$. Moreover, still outside $\Delta_d$, the diffeomorphism class of $Z(s)$ depends only on $d$. For any $s$ we equip $Z(s)$ with the restriction $g_{\omega|Z(s)}$ of the K\"ahler metric $g_\omega$. 
Then, $(Z(s), g_{\omega|Z(s)})$ can be seen as a fixed manifold with a random metric. 
In this paper we are interested in the statistics of the various metric observables of this pair, in particular the systole, the curvature, the injectivity radius and the first non-trivial eigenvalue of the Laplacian.

The \emph{systole} is the length of the smallest non contractible real curve of 
$ Z(s)$ when $n=2$, and in general the Berger $k$-systole. 
More precisely, let $(X,g)$ be a $n$-dimensional Riemannian manifold. Define 
$$\syst (X,g)= 
		\left\lbrace \begin{array}{lll}
		\inf\{\length(c), \ c \text{ non contractible} & \\ \text{smooth closed curve in } X\}  &\text{ if } \pi_1(X)\not=\{1\}\\
		 +\infty &\text{ if } \pi_1(X)=\{0\}
\end{array}
\right.
$$
and (see~\cite{berger}, \cite{gromov1983filling}) for any 
$ k\in \{1,\cdots, n \}$, 
$$ \syst_k (X,g)= 
\left\lbrace \begin{array}{ll}
\frac12 \inf\{\diam (H) | H_k(X,\R)\ni [H]\neq 0\} & \text{ if }  H_k(X,\R)\neq \{0\}\\
 +\infty &\text{ if }H_k(X,\R)=\{0\}\end{array}
\right..
$$
Our definition of the $k-$systole is different but very close to the one
given by Berger (where it is called a \emph{carcan}), which refers to the volume of the
submanifold, and not the diameter.
Note that $\syst \leq \syst_1$. 
If  $X=\C P^n$ and $E=X\times \C$, then by the Lefschetz theorem, for any odd $k\neq n-1$, any degree $d\geq 1$ and any generic section $s\in H^0(X,L^d)$, $H_k(Z(s),\R)= 0$, so that 
$\syst_{k} (Z(s))=+\infty$. This is also true for complete intersections.
The same theorem implies also that for $n\geq 3$, $Z(s)$ is simply connected, so that $\syst =+\infty$ in this case.

	 The first result of this paper provides probabilistic {lower} bounds for these systoles:                        
                        \begin{theorem}\label{sy}
                         Let $X$ be a compact smooth complex manifold of dimension $n\geq 2$ equipped with   an ample holomorphic line bundle $(L,h)\to X$  endowed with a Hermitian metric $h$ with positive curvature and K\"ahler metric $g_\omega$, and with $(E,h_E)$ a rank $r$ holomorphic vector bundle. Then, 
                        for any sequence $(a_d)_d$ converging to $0$, there exists a constant $C>0$ and a positive integer $d_0$, such that 
                        for any $k\in \{1, \cdots, 2n-2r\},$ 
                        $$ \forall d\geq d_0, \ \mathbb \mu_d\left[ \syst_k(Z(s),g_{\omega|Z(s)}) \geq 
                        \frac{a_d}{C\sqrt{\log d}} d^{-\frac{3n+1}2}
                        \right]\geq 1-C(a_d+\frac{1}d).$$
                        The same holds for the systole. Here, $\mu_d$ denotes the Gaussian measure defined by~(\ref{mesure}).
                    \end{theorem}
                In was proven~\cite[Theorem 1.16 and Corollary 1.21]{gayet2022systoles}) that there exists $c>0$ such that for any $d$ large enough,
                	\begin{itemize}
                		\item if $n-r=1$, $ \mu_d \left[\syst(Z(s),g_{\omega|Z(s)})\leq \frac{1}{\sqrt d}\right]\geq c;$
                		\item if $n$ is odd, $  \mu_d \left[\syst_{n-1}(Z(s),g_{\omega|Z(s)})\leq \frac{1}{\sqrt d}\right]\geq c.$
                	\end{itemize}

                We also give an estimate of the sectional curvatures $K$ of the random complex submanifold: 
                 \begin{theorem}\label{theoremcurvature}
                	Under the hypotheses of Theorem~\ref{sy}, for any sequence $(a_d)_d$ converging to zero,
                	there exists $C>0$ such that 
                		$$ \forall d\gg 1, \ \mathbb \mu_d\left[
                		\sup_{
                			\stackrel{
                				x\in Z(s)}
                			{ P\in \grass(2, T_xZ(s))}}|K(P)|
                		\leq \frac{C}{a_d}d^{\frac{3}2(n+1)} \log d\right]\geq  1- C(a_d+\frac{1}d).$$
                \end{theorem}
            	
                The latter and the systole estimate allow us to give a statistical control of the injectivity radius
                of the complex submanifolds:
                \begin{theorem}\label{theoremradius}
                	Under the hypotheses of Theorem~\ref{sy}, for any sequence $(a_d)_d$ converging to zero,
                	                	there exists $C>0$ such that 
	$$ \forall d\gg 1, \  \mathbb \mu_d\left[
	\inj(Z(s),g_{\omega|Z(s)}) \geq 
	\frac{a_d }{Cd^{\frac{3n+1}2}\sqrt{\log d}}
	\right]\geq  1- C(a_d+\frac{1}d).
	$$
                \end{theorem}
            Finally, we prove a control for the spectral gap. Recall that 
if $(M,g)$ is a compact smooth Riemannian manifold, if $\Delta$ its
            associated Laplace-Beltrami operator, then
            its spectrum is positive, infinite, discrete and unbounded.
            The first non vanishing eigenvalue $\lambda_1= \lambda_1(M,g)$ is called the \emph{spectral gap} and is of great importance in various problems. For instance, the heat kernel decreases with the speed $\sqrt{\lambda_1}.$
            
       \begin{theorem}\label{theoremeigenvalue}
      	Under the hypotheses of Theorem~\ref{sy}, for any sequence $(a_d)_d$ converging to zero,
      	there exists $C>0$ such that 
     	$$ \mathbb \mu_d\left[
      	\lambda_1(Z(s),g_{\omega|Z(s)})
      	\geq \exp(-\frac{C}{\sqrt{a_d}} d^{\frac{3n+3+12r}4}\sqrt{\log d})\right]\geq  1- C(a_d+\frac{1}d).$$
      \end{theorem}
In the proof of Theorem~\ref{theoremeigenvalue}, we use a new deterministic bound 
for the diameter of a degree $d$ complex submanifold, whose proof is close to the one in~\cite{feng1999diameter}. For plane curves, our bound on the diameter is $Cd^3$, while in ~\cite{feng1999diameter} it is $Cd^4$.
\begin{theorem}\label{theorem: diameter} Under the hypotheses of Theorem~\ref{sy}, there exists $c=c(X,L,\omega)$ such that  
	$$\forall d\geq 1, \ 
\forall s\in H^0(X,E\otimes L^d), \ 
	\mathrm{diam}(Z(s),g_{\omega |Z(s)})\leq cd^{3r}.$$
\end{theorem}

\noindent
{\bf Open questions:} what are the laws of the various metric observables, in particular $\syst(Z(s))$ and $\diam (Z(s))$? The lower bound for the spectral gap is pretty bad compared to what is known about hyperbolic surfaces, can it be amended, at least with a polynomial bound?  \\

\noindent
{\bf Idea of the proofs for $n=2$ and for the standard setting}. The systole estimate holds on the following simple idea: if for any point $x\in Z(P)$, the complex curve $Z(P)$ is a locally a graph in $\C P^2$ over a complex disc of size $r$ centered at $x$, then any non-contractible curve in $X$ has diameter larger than $r$. For any $x\in X$, a quantitative implicit function provides a lower bound for $r$ in terms of an upper bound of $\mathscr{C}^2$ norm of $P$ and a lower bound of its gradient over $Z(P)$, see Corollary~\ref{coro} and Proposition~\ref{loco}. Now, this lower bound is controlled by the distance of $P$ to the discriminant locus $\Delta_d$, see Proposition~\ref{burg}. Moreover, one can estimate the probability that $P$ is far enough to $\Delta_d$, see Proposition~\ref{burg}. The fact that $\Delta_d$ is an algebraic set is important here. On the other hand, crude estimates for the $\mathscr{C}^2$ norm can be achieved by deterministic arguments, see Lemma~\ref{deter}, and far better ones in probability by L\'evy concentration arguments, see~Theorem~\ref{borne}. The curvature estimates are given also by a control of the gradient and the $\mathscr{C}^2$ norm, see Proposition~\ref{curv}. The injectivity radius estimate is a direct consequence of the two former bounds, see Theorem~\ref{ehrlich}. Finally, the spectral estimate is a consequence of a theorem by Gromov, see Theorem~\ref{eig}, of the estimate of the curvature, and of a deterministic upper bound for the diameter, see Theorem~\ref{dia}.\\

\noindent
{\bf Acknowledgments.}
The research leading to these results has received funding from the French Agence nationale de la ANR-20-CE40-0017 (Adyct).
\section{Preliminaries}

\subsection{Asymptotics of the Bergman kernel}\label{bergman}
In this paragraph we assume that the setting and hypotheses of Theorem~\ref{sy} are satisfied. The covariance function $E_d$ for the Gaussian field generated by the holomorphic sections $s\in H^0(X,E\otimes L^d)$ is defined by
$$ \forall z,w\in X, 
\ E_d (z,w) = \mathbb E \left[s(z)\otimes (s(w))^*\right]\in (E\otimes L^d)_z\otimes (E\otimes L^d)_w^*,$$
where the averaging is made for the measure $\mu_d$ given by~(\ref{mesure}), where $E^*$ is the (complex) dual of $E$ and 
$$\forall w\in X, \ \forall s,t\in (E\otimes L^d)_w, \  s^* (t) = \langle s,t\rangle_{h_d(w)}.$$
The covariance $E_d$ is the \emph{Bergman kernel}, that is the kernel of the orthogonal projector from $L^2(M,E\otimes L^d)$ onto $H^0(M,E\otimes L^d)$. This fact can be seen through the equations
$$\forall z,w\in M, \ E_d(z,w) = \sum_{i=1}^{N_d} S_i (z)\otimes S_i^*(w),$$
where $(S_i)_i$ is an orthonormal basis of $H^0(M,E\otimes L^d)$ for the Hermitian product~(\ref{prod}).
Recall that the metric $g_\omega$ is induced by the curvature form $\omega$  and the complex structure. 
It is now classical  that the Bergman kernel has a universal rescaled (at scale $\frac{1}{\sqrt d}$) limit, the Bargmann-Fock kernel $\mathcal P$:
\begin{equation}\label{mp}
\forall z,w\in \C^n, \ \mathcal P(z,w):=  \exp\left(
-\frac{\pi}2 (\|z\|^2 +\|w\|^2 -2\langle z,w\rangle)
\right).
\end{equation}
 Theorem~\ref{Dai} below quantifies this phenomenon. For this, we need to introduce local trivializations and charts. 
Let $x\in X$ and $R>0$ such that $2R$ is less than the radius of injectivity of $X$ at $x$. Then the exponential map based at $x$ induces a chart near $x$ with values in $B_{T_xX}(0,2R)$. We identify a point in $X$ with its coordinates. The parallel transport provides a trivialization 
$$\varphi_x : B_{T_xX}(0,2R)\times (E\otimes L^d)_{x}\to (E\otimes L^d)_{|B_{T_xX}(0,2R) }$$ 
which induces a trivialization of 
$(E\otimes L^d)\boxtimes
(E\otimes L^d)^*_{|B_{T_xX}(0,2R)^2}$.
Under this trivialization, the Bergman kernel 
$ E_d$ becomes a map from $T_x M^2 $ with values into $\End\left((E\otimes L^d)_x\right)$.

\begin{theorem}(\cite[Theorem 1]{ma2013remark})\label{Dai} Under the hypotheses of Theorem~\ref{sy}, let $m\in \Nn$. Then, there exist $C>0$, such that for any $k\in\{0,\cdots, m\}, $ for any $x\in X$, 
	$\forall z,w\in B_{T_xX}(0,\frac{1}{\sqrt d}), $ 
	\beqe \left\|	D^k_{(z,w)}\left(
	\frac{1}{d^n}E_d(z,w) - \mathcal P(z\sqrt d,w\sqrt d )
	\ \id_{(E\otimes L^d)_x}
	\right)
	\right\|	
	\leq C d^{\frac{k}2-1}.
	\eeqe
\end{theorem}
The original reference is more general, see~\cite[Proposition 3.4]{letpuch} for the present simplification.

\subsection{Upper bounds for $ \mathscr{C}^k$-norms}

In this paragraph we provide 
probabilistic estimates for the $ \mathscr{C}^k$ norm. 
We begin by the following crude deterministic estimate:
\begin{lemma}\label{deter} Under the hypotheses of Theorem~\ref{sy}, for all $k\in \{0,1,2\}$, there exists $C>0$ such that 
	$$ \forall d\gg 1, \ \forall s\in  S^d, \ \sup_{x\in X} \|\nabla^k s(x)\|\leq Cd^{\frac{n+k}2}.$$
\end{lemma}
\begin{preuve} We prove this estimate for $k=1$, the other cases are similar. 
	Let  $(S_i)_{i=1\cdots N_d}$ be an orthonormal basis  of $ H^0(X,E\otimes L^d)$ and $s=\sum_{i=1}^{N_d}a_i S_i\in S^d$, where $(a_i)_i \in \C^{N_d}\cap \mathbb S^{2N_d-1}.$ By Cauchy-Schwarz, 
	\begin{eqnarray*} 
		\forall x\in X, \ 
		\|\nabla s(x)\|^2 &\leq &(\sum_{i=1}^{N_d} a_i \|\nabla S_i(x)\|)^2 \leq  \sum_{i=1}^{N_d} \|\nabla S_i(x)\|^2.
	\end{eqnarray*}
	The latter equals
	$  \| \sum_{i=1}^{N_d} \nabla S_i(z)\otimes \nabla^* S_i (w)_{|z=w=x}\|$, which is equal to $\|\nabla_z \nabla^*_w E_d\|^2_{z=w=x}.$
	Theorem~\ref{Dai} now concludes. 
\end{preuve}
In probability, we can get better estimates.
\begin{theorem}\label{borne} Under the hypotheses of Theorem~\ref{sy}, for any $p \in \Nn$ and $k\in \{0,1,2\}$, 
	there exists $C>0$, such that 
	$$\forall d\gg 1, \ \mu_d \left \lbrace
	s\in  S^d, \ \sup_{x\in X} \|\nabla^k s(x)\|> C\sqrt{d^k \log d}
	\right\rbrace
	\leq \frac{1}{d^p},$$
where $S^d = \mathbb S H^0(X,E\otimes L^d)$. 
\end{theorem}
This theorem has been proven by Shiffman and Zelditch in the case of a line bundle, see~\cite[Theorems 1.1 and 1.4]{SZ}. For the sake of completeness, we give a proof for our more general setting.
\begin{preuve}[ of Theorem~\ref{borne}]
	Let  $(S_i)_{i=1\cdots N_d}$ be an orthonormal basis  of $ H^0(X,E\otimes L^d)$.
Fix $x\in X$, and let $(h_\ell(x))_{\ell\in \{1, \cdots, r\} }$ be an orthonormal basis of $E\otimes L^d_{|x}$. For $\ell\in \{1, \cdots, r\}$, define
$$
\sigma_\ell(x)= (
\sigma^i_\ell(x)
)_{i\in \{1, \cdots, N_d\}}:=
\left(\langle S_i(x),h_\ell(x)\rangle_{h_d}\right)_{i\in \{1, \cdots, N_d\}}
.$$
	Since $L$ is ample, there is no base point for $H^0(X,E\otimes L^d)$ for $d$ large enough. Consequently, $$(\sigma_{\ell}(x))_{\ell\in \{1, \cdots, r\}}\neq 0.$$
Let $a:=(a_i)_{i\in \{1, \cdots, N_d\}}\in \C^{N_d}$ and  $s=\sum_{i=1}^{N_d}a_i S_i\in H^0(X,E\otimes L^d)$. Then, 
	$$  |s(x)|_{h_{d}}^2= \sum_{\ell=1}^r
	\left|\langle a,
\sigma_\ell(x)\rangle_{\C^{N^d}}\right|^2$$
  By the L\'evy concentration, see~\cite[Lemma 2.1 and (23)]{SZ}, for any $\ell\in \{1, \cdots, r\}$
  such that $\sigma_\ell(x)\neq 0$,  for any $C>0$, 
	$$\mu_d \left \lbrace
	s\in  S^d, \  
	\frac{|\langle a, \sigma_{\ell}(x)\rangle|}{|\sigma_\ell(x)|}> {C}d^{-n/2}\sqrt{\log d}
	\right\rbrace\leq
	\exp\left(-(N_d-1){C^2}\frac{\log d}{d^n}\right).
	$$
Hence,
\begin{equation}\label{kasos} \mu_d \left \lbrace
	s\in  S^d, \  \frac{|\langle a, \sigma_{\ell}(x)\rangle|}{(\sum_{\ell'=1}^r |\sigma_{\ell'}(x)|^2)^{1/2}}> {C}d^{-n/2}\sqrt{\log d}
	\right\rbrace\leq 
	\exp\left(-(N_d-1){C^2}\frac{\log d}{d^n}\right),
	\end{equation}
	and the inequality still holds if $\sigma_\ell(x)=0$. 
Since for any $\alpha>0$,
$$  \{s\in S^d, \ |s(x)|_{h_d}>\alpha\} \subset \bigcup_{\ell\in \{1,\cdots, r\}}\{s\in S^d, \ |\langle a, \sigma_{\ell}(x)\rangle|> \alpha/\sqrt r\},$$
%
(\ref{kasos}) implies that 
	$$ \mu_d \left \lbrace
	s\in  S^d, \  \frac{|s(x)|}{(\sum_{\ell=1}^r |\sigma_\ell(x)|^2)^{1/2}}> Cd^{-n/2}\sqrt{\log d}
	\right\rbrace\leq r
	\exp\left(-(N_d-1)\frac{C^2}{r}\frac{\log d}{d^n}\right).$$
By Theorem~\ref{Dai}, for any $x\in X$,
 $$\sum_{\ell=1}^r |\sigma_\ell(x)|^2= \sum_i |S_i(x)|_{h_d}^2 = E_d(x,x)\sim_d d^n,$$
 hence for any $C>0$, for $d$ large enough and any $x\in X$, 
 $$ \mu_d \left \lbrace
 s\in  S^d, \  |s(x)| > C\sqrt{\log d}
 \right\rbrace\leq r
 \exp\left(-(N_d-1)\frac{C^2}{2r^2}\frac{\log d}{d^n}\right).$$
	Now, since the connection $\nabla$ is metric,
$$ \forall x\in X,\  D|s(x)|^2_{h_d}=
	2\Re \langle \nabla s(x), s(x)\rangle_{h_d}.
$$
	By Lemma~\ref{deter} this implies that
	there exists $c>0$, such that for any $d$, 
		\begin{equation}\label{bask}
		 \|D(|s|_{h_d})\|_\infty \leq cd^{{\frac{1+n}2}}.
		\end{equation}
	As in~\cite{SZ}, we cover $X$ by a collection of $k_d$ balls $(B(x_j,r_d))_{1\leq j\leq k_d}$ where
	$ r_d = d^{-\frac{n+1}2}$. We can choose the centers optimally such that 
	there exists $c>0$ independent of $d$ such that
	$$ k_d \leq  cd^{	n(n+1)}.$$
	By~(\ref{bask}) and the mean value inequality, we obtain 
	$$ \forall C>0, \ \forall d\gg 1, \ \mu_d \left \lbrace
	s\in  S^d, \  \|s\|_{\infty}> C \sqrt{\log d}
	\right\rbrace\leq cd^{n(n+1)}
	\exp\left(-(N_d-1)\frac{C^2}{2r^2}\frac{\log d}{d^n}\right).$$
	Since there exists $c>0$, such that $N_d \sim c d^n$, choosing $C$ large enough gives
	the result for $k=0$. The cases $k=1,2$ are proven similarly, see~\cite{SZ}.
\end{preuve}

\subsection{Discriminant and quantitative transversality} In this paragraph we assume that the setting and hypotheses of Theorem~\ref{sy} are satisfied. For any  positive integer $d$ and 
any $s\in H^0(X,E\otimes L^d)$, define 
$$ \dist_{\Delta_d}(s):= \min_{\tau\in \Delta_d}\|s-\tau\|_{L^2},$$ 
where $\Delta_d$ denotes the discriminant locus of $H^0(X,E\otimes L^d)$, see the introduction.
For any positive $d$, let $\nabla$ be the Chern connection associated to $(E\otimes L^d,h_d)$.

\begin{proposition}\label{raf2}
	Under the hypotheses of Theorem~\ref{sy},
	for any  positive integer $d$ and 
	any $s\in H^0(X,E\otimes L^d)$,
	$$ \dist_{\Delta_d} (s) =
	\min_{x\in X} \left(
	\frac{|s(x)|^2_{h_d}}{d^n}+\min_{M\in \mathcal L_{sing}(T_xX, (E\otimes L^d)_x)}
	\frac{\|\nabla s(x)-M\|^2_{h_d}}{d^{n+1}}
	\right)^{1/2}
	(\pi^{n/2} +O(\frac{1}{d})),$$
	where the error term is independent of  $s$, and $\mathcal L_{sing}$ denotes the space of morphisms which are not onto.
\end{proposition}
\begin{preuve}
The proof of this proposition follows the lines of the proof of~\cite[Lemma 3.8]{ancona-complete} where the case $E=X\times \C^r$ was treated. 
\end{preuve}

We need to bound the degree of the discriminant locus $\Delta_d$.
	\begin{proposition}\label{dede}
		Under the hypotheses of Theorem~\ref{sy}, the discriminant locus $\Delta_d\subset  H^0 (X,E\otimes L^d)$ is an algebraic hypersurface and its degree satisfies 
		$$ \deg(\Delta_d)\underset{d\to \infty}\sim r\binom{n+1}{r}d^n.$$ 
		\end{proposition}
	\begin{preuve} For $E=X\times\C^r$, this Proposition is \cite[Lemma 2.3]{ancona-complete}. We give the proof for the general case for the sake of completeness.
	By hypothesis, $L$ is ample and then  $E\otimes L^d$ is $1$--jet spanned for $d$ large enough. Moreover, the incidence variety $$\Sigma=\{(x,[s])\in X\times \mathbb{P}(H^0(X,E\otimes L^d)) | \hspace{1mm}x \mathrm{\hspace{1mm}is\hspace{1mm} a\hspace{1mm} singular\hspace{1mm} point \hspace{1mm}of\hspace{1mm}} Z(s)\}$$ is also smooth as soon as $d$ is large enough. Remark that the image of $\Sigma$ in $\mathbb{P}(H^0(X,E\otimes L^d))$ under the second projection is precisely the discriminant locus $\Delta_d$. Moreover, the projection $\Sigma\rightarrow \Delta_d$ is a  birational map (that is, a generic section $s$ of the discriminant only has one singular point). We are then in the hypothesis of \cite[Proposition 1.1 and Corollary 2.4]{Abo_Lazarsfeld_Smith_2022} that gives us the following equality for the degree of $\Delta_d$:
		\begin{equation}\label{segre} \deg \Delta_d = \int_X p_*\left(
	\{ c(p^* \Omega_X) s (p^* (E\otimes L^d)^*\otimes  \mathcal{O}_{\mathbb{P}(M^*)}(-1))
	\}_{n-r+1}c_1(\mathcal{O}_{\mathbb{P}(M^*)}(1))^{N_d-2}\right).
	\end{equation}
		Here we have used the following notation: 	 $\Omega_X$ denotes the complex cotangent bundle of $X$, $$p : X\times \mathbb{P}(H^0(X,E\otimes L^d)^*)\to X$$
	denotes the projection onto the first factor, $M$ is the fiber bundle over $X$ with fiber 
	$M_x = \{s\in H^0, s(x)=0 \},$ $c(\cdot)$ is the total Chern class, $s(\cdot)$ is the total Segre class, $N_d$ is the dimension of $H^0(X,E\otimes L^d)$ and $\{\cdot\}_\ell$ denotes the degree $\ell$ part of a non--pure cohomology class.
	
	We now develop the degree $n-r+1$ part of $s(p^* (E\otimes L^d)^*\otimes  \mathcal{O}_{\mathbb{P}(M^*)}(-1))$. Let us denote by $F$ the rank $k$ complex vector bundle $E\otimes L^d$.	
	The total Segre class of $p^*F^*\otimes\mathcal{O}_{\mathbb{P}(M^*)}(-1)$ can be decomposed as $$s(p^*F^*\otimes\mathcal{O}_{\mathbb{P}(M^*)}(-1))=\sum_{i=0}^rs_i(p^*F^*\otimes\mathcal{O}_{\mathbb{P}(M^*)}(-1)).$$ Moreover, each Segre class $s_i(p^*F^*\otimes\mathcal{O}_{\mathbb{P}(M^*)}(-1))$	can be itself be decomposed as 
	$$s_i(p^*F^*\otimes\mathcal{O}_{\mathbb{P}(M^*)}(-1))= \sum_{j=0}^i (-1)^{i-j} \binom{r+i-1}{r+j-1} s_j(p^*F^*) c_1(\mathcal{O}_{\mathbb{P}(M^*)}(-1))^{i-j}$$
	see for example~\cite[Example 3.1.1]{fulton}.
	We then obtain that the degree $n-r+1$ part of $c(p^* \Omega_X) s (p^* (E^*\otimes L^{-d})\otimes \mathcal{O}_{\mathbb{P}(M^*)}(-1))$ equals
	\begin{multline} \sum_{k+h=n-r+1}c_k(p^* \Omega_X)s_h(p^* F^*\otimes \mathcal{O}_{\mathbb{P}(M^*)}(-1)) \\
	= \sum_{k+h=n-r+1}c_k(p^* \Omega_X)\sum_{j=0}^{h}(-1)^{h-j} \binom{r+h-1}{r+j-1} s_j(p^*F^*) c_1(\mathcal{O}_{\mathbb{P}(M^*)}(-1))^{h-j}\\
	= \sum_{k+h=n-r+1}c_k(p^* \Omega_X)\sum_{j=0}^{h} \binom{r+h-1}{r+j-1} s_j(p^*F^*) c_1(\mathcal{O}_{\mathbb{P}(M^*)}(1))^{h-j}.
	\end{multline}
	This implies that the degree of the discriminant equals 
	\begin{multline} \deg(\Delta_d)=\int_Xp_*\left( \sum_{k+h=n-r+1}c_k(p^* \Omega_X)\sum_{j=0}^{h}\binom{r+h-1}{r+j-1} s_j(p^*F^*) c_1(\mathcal{O}_{\mathbb{P}(M^*)}(1))^{h-j+N_d-2}\right) \\
	=\displaystyle\int_X\sum_{k+h=n-r+1} c_k(\Omega_X)\sum_{j=0}^{h}\binom{r+h-1}{r+j-1} s_j(F^*) p_*c_1(\mathcal{O}_{\mathbb{P}(M^*)}(1))^{h-j+N_d-2}\\
	=\displaystyle\int_X\sum_{k+h=n-r+1} c_k(\Omega_X)\sum_{j=0}^{h}\binom{r+h-1}{r+j-1} s_j(F^*)c_{h-j-1+r}(F),
	\end{multline}
	where in the last equality we have used that $p_*c_1(\mathcal{O}_{\mathbb{P}(M^*)}(1))^{h-j+N_d-2}=c_{h-j-1+r}(F)$ (see for example \cite[Page 4]{Abo_Lazarsfeld_Smith_2022}).
	For $F=E\otimes L^d$, the last integral equals 
	\begin{multline} \displaystyle\int_X\sum_{k+h=n-r+1} c_k(\Omega_X)\sum_{j=0}^{h}\binom{r+h-1}{r+j-1} s_j(E^*\otimes L^{-d})c_{h-j-1+r}(E\otimes L^{d})
\\
=\displaystyle\int_X\sum_{k+h=n-r+1} c_k(\Omega_X)\sum_{j=0}^{h}\binom{r+h-1}{r+j-1} \times \\
\times \sum_{i=0}^j(-1)^{j-i}\binom{r+j-1}{r+i-1}s_i(E^*)c_1 (L^{-d})^{j-i}\sum_{m=0}^{h-j+1+r}\binom{j+m+1-h}{m}c_{h-j-1+r-m}(E)c_1(L^d)^m.
	\end{multline}
	where in the last equality we have used that $c_{\ell}(E\otimes L^{d})=\sum_{m=0}^{\ell}\binom{r-\ell+m}{m}c_{\ell-m}(E)c_1(L^d)^m$, see for example~\cite[Example 3.2.2]{fulton}.
Using that $c_1(L^d)=dc_1(L)$, we can see the previous integral as the following polynomial in $d$ whose coefficient are integral of Segre and Chern classes of $\Omega_X, L$ and $E$:
\begin{multline}
\displaystyle\int_X\sum_{k+h=n-r+1} c_k(\Omega_X)\sum_{j=0}^{h}\binom{r+h-1}{r+j-1} \times \\
\times \sum_{i=0}^j\binom{r+j-1}{r+i-1}s_i(E^*)c_1 (L)^{j-i}\sum_{m=0}^{h-j-1+r}d^{m+j-i}\binom{j+m+1-h}{m}c_{h-j-1+r-m}(E)c_1(L)^m.
	\end{multline}

 We can see that the degree of this polynomial is $n$, and the degree $n$ coefficient of this polynomial can be obtained by  looking at the indices $h=n-r+1, k=0, m+j=n$ and $i=0$. Such coefficient is then equal to $\int_X\sum_{j=0}^{n-r+1}\binom{n}{r+j-1}\binom{r+j-1}{r-1}\binom{r}{n-j}c_1(L)^n$.
 In the sum, only the terms associated to $j=n-r$ and $j=n-r+1$ are non-zero, so that the coefficient equals, after some simplifications to $r\binom{n+1}{r}.$
\end{preuve}
The following proposition is a generalization of~\cite[Lemma 3.4]{ancona-exp} in our setting:
\begin{proposition}\label{burg} Under the hypotheses of Theorem~\ref{sy}, there exists $c, C>0$ and a positive integer $d_0$, such that for any sequence $(r_d)_{d\in \Nn^*}$
	satisfying $\forall d\geq d_0, \ r_d\leq cd^{-2n}$, 
	$$ \forall d\geq d_0, \ \mu_d \left\lbrace
	s\in  S^d, \ \dist_{\Delta_d} (s)
	\leq r_d \right\rbrace\leq Cr_dd^{2n},$$ 
	where $S^d\subset H^0(X,E\otimes L^d)$ denotes the $L^2$-unit sphere and $\mu_d$ is the measure defined by~(\ref{mesure}).
\end{proposition}
\begin{preuve}[ of Proposition~\ref{burg}]
	This proposition was proved in the real setting, but its proof holds in the complex setting, \emph{mutatis mutandis}. The main tool is~\cite[Theorem 21.1]{burgisser}, which estimates this volume in terms of the dimension, growing like $O(d^n)$, and the degree of the discriminant locus $\Delta_d$. By Proposition~\ref{dede}, this degree grows like $O(d^n)$.
\end{preuve}

\subsection{Local graphs}
For the systole and the injectivity radius estimates we will use a quantitative version of the implicit function theorem:
\begin{proposition}\label{tfi}
	Let $p,q$ be positive integers, $f : \R^{p}\times \R^q \to \R^q$ be a $\mathscr{C}^2$ function and $x = (x_1, x_2 ) \in  \R^{p}\times  \R^q$ such
	that $f (x) = 0$ and the partial derivative $D_{x_2} f (x) : \R^q \to \R^q$ is invertible. Choose $\delta>0$ such that
	$$\sup_{y\in B(x,\delta)}
	\|\id_{|\R^{q}} - (D_{x_2} f (x))^{-1}
	D_{x_2} f (y) \|\leq  1/2.$$
	Let $\displaystyle C = \sup_{y\in B(x,\delta)} \|D_{x_1} f (y)\|,$ $  M = \|(D_{x_2} f (x))^{-1}\|, $  
	$\delta' = \delta(2MC)^{-1}$
	and $$I_{\delta'} := \{y_1 \in \R^{p} : |y_1 - x_1 |<\delta'\}.$$ Then there exists $\varphi: I_{\delta'} \to \R^p$  a $\mathscr{C}^2$ function such that
	$$ \forall y = (y_1 , y_2 ) \in \R^{p}\times \R^q,
	|y_1 - x_1| < \delta',\ 
	|y_2 - x_2 | < \delta,\ 
	f (y) = 0 \Leftrightarrow y_2 = \varphi(y_1 ).$$
\end{proposition}
The proof of this proposition is the classical one, keeping track the parameters. 
We will use the following corollary:
\begin{corollary}\label{coro} Let $1\leq r\leq n$ be integers,  $f : 2\mathbb B \subset \C^n \to \C^r$ be a $\mathscr{C}^2$ function. 
	Let $\epsilon, C_1, C_2$ be positive constants such that $\|Df\|_{\infty}\leq  C_1$ and 
$	\|D^2f\|_{\infty}\leq  C_2.$
	Let  $x \in B$ be satisfying 
	$f (x) = 0$ and $$
	\min_{M\in \mathcal L_{sing}(\C^n, \C^r)}\|Df(x)-M\| > \epsilon.$$
	Then there exist complex orthonormal coordinates $z=(z_1, \cdots z_n)$ such that if $z'=(z_1, \cdots, z_{n-r})$, $z''=(z_{n-r+1}, \cdots, z_n)$ and 
	$$\delta = \frac{\epsilon^2}{4C_1C_2},$$
	the zero set $Z(f)\cap \left( B(x', \frac{\epsilon}{2C_2})\times B(x'', \delta)\right)$ is the graph over $B(x', \frac{\epsilon}{2C_2}) $ of a $\mathscr{C}^2$ function, where $x=(x',x'') $ in the new coordinates.
\end{corollary}

\begin{lemma}\label{patate}Under the hypotheses of Theorem~\ref{sy}, there exists $C>0$ such that for any $d$ large enough and any $x\in X$, 
	under the trivialization above
	for any local $\mathscr{C}^2$ section $s$ of $E\otimes L^d$,
	\begin{eqnarray*} \|Ds - \nabla s\|_{|y}& \leq & C d|y||s(y)|_{h_d}\\
		\|D^2 s -\nabla^2 s\|_{|y}&\leq &C\left((d+d^2 |y|^2) |s|_{h_d}+
		d|y|\|\nabla s\|\right)
	\end{eqnarray*}
\end{lemma}
\begin{preuve}If $e$ (resp. $(\epsilon_\ell)_{\ell\in \{1, \cdots, r\}}$) is a trivialization of $L$ (resp.  of $E$), $A_L$ (resp. $A_E$) the associated connection 1-form, then for any $d\geq 1$ and $\mathscr{C}^2$ local function $f$, for any $\ell \in \{1,\cdots, r\}$,
	\begin{equation}\label{piaf} \nabla ( f e^d\otimes \epsilon_\ell) = 
	(Df+dfA_L + fA_E) e^d\otimes \epsilon_\ell,
	\end{equation}
	see also~\cite[(1.6.21)]{ma2007holomorphic}. Since the trivialization is done by radial parallel transport $A_E=O(|y|)$ and $A_L=0(|y|)$, hence the first assertion. The second assertion is proven by differentiating~(\ref{piaf}).
\end{preuve}

The following probabilistic  Lemma provides
a global scale at which 
the submanifold can be seen everywhere
at as a graph with controlled parameters.

\begin{lemma}\label{loco}
	Under the hypotheses of Theorem~\ref{sy}, let $(a_d)_d$ be a sequence converging to 0. Then, there exists $C>0$ and a positive integer $d_0$ such that for any $d\geq d_0$, if 
	$$ 
 \rho_d = \frac{1}{C\sqrt{\log d}} a_d d^{-\frac{3n+1}2}, 
	$$ 	with probability at least 
		$ 1 - C(a_d+\frac{1}d), $ a section $s\in  S^d$ satisfies the following event: 
	for any $x\in Z(s)$, under the local chart and trivialization defined in \S~\ref{bergman}, 
	\begin{eqnarray*}
\forall y\in B(x,\rho_d), \ |Ds(y)|&\leq &
	C	\sqrt{d\log d}\\
	\min_{M\in \mathcal L_{sing}(T_xX, E\otimes L^d_{|x})}\|Ds(x)-M\|&\geq& 		\frac{1}{2} a_d d^{\frac{-3n+1}2} \\
	|D^2s(y)|&\leq &
	Cd	 \sqrt{\log d}.
	\end{eqnarray*} 
\end{lemma}
\begin{preuve}
Let $(a_d)_d$ be a sequence converging to 0, $A$ be the event 
	$$ A =\{s\in S^d, \ 
\dist_{\Delta_d} (s)
\geq a_d d^{-2n}\}$$
and for $k=\{0,1,2\}$, $B_k$ be the event
$$ B_k = \left\{s\in S^d, \ C
\sup_X \|\nabla^k s\|\leq  C\sqrt{d^k \log d}
\right\}$$
and $B =  A\cap B_0 \cap B_1 \cap B_2.$
	By Proposition~\ref{burg} and Theorem~\ref{borne}, 
there exists $C$ depending on the sequence $(a_d)_d$ such that 
$$\forall d\geq 1, \  \mu_d (B) \geq 1- C(a_d+\frac{1}d).$$ 
From now on, assume that $s\in B$. 
By Proposition~\ref{raf2}, since  $s\in A$, 
there exists  a universal positive constant $C'$ such that 
\begin{equation}\label{deadman}\min_{M\in \mathcal L_{sing}(T_xX, (E\otimes L^d)_x)}\|\nabla s(x)-M\|^2_{h^d}\geq (a_d d^{-2n})^2 d^{n+1}
- d|s(x)|_{h_d}^2.
\end{equation}
In the trivialization, Lemma~\ref{patate} implies that
\begin{eqnarray*}
 \forall y\in B(x,\frac{1}{\sqrt d}),\ \| Ds\|&\leq &\|\nabla s\| + d|y| |s|_{h_d}\leq C\sqrt {d\log d}\\
 \|D^2 s\|&\leq &Cd\sqrt{\log d}.
 \end{eqnarray*}
By the mean value inequality and $x\in Z(s)$, the estimate above provides the estimate:
\begin{equation}\label{godzilla} \forall y \in B(x,\frac{1}{\sqrt d}), \
|s(y)|_{h_d}\leq C|y|\sqrt{d\log d}.
\end{equation}
Lemma~\ref{patate}, (\ref{deadman}) and (\ref{godzilla})  imply that
\begin{eqnarray*}
	\forall y \in B(x,\frac{1}{\sqrt d}), \ 
 \min_{M\in \mathcal L_{sing}(T_xX, E\otimes L^d_{|x})}\|Ds(x)-M\|
&\geq &
a_d d^{\frac{-3n+1}2}- C(\sqrt d + d|y|)|y| \sqrt{d\log d} \\
&\geq & a_d d^{\frac{-3n+1}2}- Cd|y| \sqrt{\log d}
\end{eqnarray*} 
Hence, for $$y\in B\left(x,\frac{1}{2C\sqrt{\log d}} a_d d^{-\frac{3n+1}2}\right)$$
 we have for $d$ large enough independent of $y$ and $x$,
\begin{eqnarray*}
\min_{M\in \mathcal L_{sing}(T_xX, E\otimes L^d_{|x})}\|Ds(x)-M\|&\geq& 		\frac{1}{2} a_d d^{\frac{-3n+1}2}.
\end{eqnarray*}
\end{preuve}

\section{Proofs of the theorems}
 \subsection{Systolic estimates}
 \begin{preuve}[ of Theorem~\ref{sy}] 
 	We fix a sequence $(a_d)_d$ which converges to 0. Then, there exist a positive constant $C$ and a positive integer $d_0$ such that for any $d\geq d_0$, 
 	by Lemma~\ref{loco} and Corollary~\ref{coro},
with probability at least 
$$ 1 - C(a_d+1/d), $$ 
for any $x\in X$, there exists local complex coordinates $z=(z',z'')\in \C^{n-r}\times \C^r$ isometric at $x$, such that 
on $$
B\left(x, \frac{a_d}{C\sqrt{\log d}} d^{-\frac{3n+1}2}\right)
\bigcap \left[
B\left(x',\frac{a_d}{4C\sqrt {\log d} } d^{-\frac{3n+1}2}\right)
\times
B\left(x'',\frac{a^2_d}{16C^2 \log d }  d^{{-3n+\frac12}}\right)
\right],$$
 $Z(s)$ is a graph of a complex function over the second ball with values in the third ball.
 The radii of the two first balls above have are equal up to a constant, and the smallest between the first and the third radius is the third.  Hence, any topologically non trivial submanifold passing through $x$ has a diameter larger than the first radius, hence the result.
	\end{preuve}

\subsection{Sectional curvatures}

Let $M$ be a submanifold of the Riemannian manifold $(N,g)$, $x\in M$ and let
\begin{eqnarray*}
	\sigma : T_x M\times T_x M &\to& N_x M\\
	(X,Y)& \mapsto & (\nabla_X Y)^\perp
\end{eqnarray*}
where $NM\subset TN$ denotes the normal bundle over $M$ and $\nabla$  the Levi-Civita connection associated to $g$. 
\begin{proposition}\label{gauss}(Gauss's equations~\cite[Theorem 3.6.2]{jost2008riemannian})
Let $M$ be a submanifold of dimension $m$ of the Riemannian $p$-dimensional manifold $N$, $x$ a point of $ M$.  
Then,
\begin{eqnarray*}
 \ \forall X,Y,Z,W\in T_xM,\ 
	\langle R^M(X,Y)Z,W\rangle&=& \langle R^N(X,Y)Z,W\rangle  +
	\\ &&\langle \sigma(Y,Z),\sigma(X,W)\rangle- \langle \sigma(X,Z),\sigma(Y,W)\rangle,
\end{eqnarray*}
where $R^M$ and $R^N$ denote the Riemannian curvature tensor for $(M,g_{|M})$ and $(N,g)$ respectively. 
\end{proposition}
Recall that for any $x\in M$ and any plane $P\subset T_xM$ spanned by an orthonormal basis $\{X,Y\}$, 
the sectional curvature of $P$ is defined by 
\begin{equation}\label{cas}
K(P) = \langle R^M(X,Y)Y,X\rangle.
\end{equation}
	The following proposition computes the curvature of the vanishing locus of a transverse section. It is very general, quite simple but we could not find any reference where it is written. Hence, it has its own interest. 
\begin{proposition}\label{Marilyn}Let $1\leq r\leq n$, $(M,g)$ be a smooth Riemannian manifold of dimension $n$, $E\to M$ be a rank $r$ smooth real vector bundle equipped with a metric $h$ and a metric connection $\nabla$. Let $s\in C^\infty(M,E)$ be a smooth section vanishing transversally on $Z(s)\subset M$. Then,
	for any $x\in Z(s)$,
	\begin{eqnarray*}
		\ \forall X,Y,Z,W\in T_xZ(s),\ 
		\langle R^{Z(s)}(X,Y)Z,W\rangle&=& \langle R^M(X,Y)Z,W\rangle + \\
&&		( \nabla s G \nabla s^*)^{-1}(\nabla^2_{Y,Z}s) (\nabla^2_{X,W} s)- \\&&
		( \nabla s G \nabla s^*)^{-1}(\nabla^2_{X,Z}s) (\nabla^2_{Y,W} s),				
	\end{eqnarray*}
where all the arguments are computed at $x$ and $\nabla s^*=(\nabla s(x))^*\in L(E^*,T_xM^*)$ denotes the adjoint of $\nabla s(x)$, and 
where $\nabla^2_{VW} = \nabla_V\nabla_W-\nabla_{\nabla_V W}$
denotes the second covariant derivative of $\nabla$.
\end{proposition}
Recall that the definition of $\nabla^2$ involves the connection $\nabla$ and the Levi-Civita connection associated to $g$.
\begin{preuve}
	Let $x\in X$ and 
	$(t_i)_{i\in \{1, \cdots, r\}}$ be a local orthonormal (for $h$) 
	frame of $E$, $(z_i)_{i\in \{1, \cdots n\}}$ be local holomorphic coordinates, such that $(\partial_{ z_i})_{i=1,\cdots, n}$ are orthonormal at $x$ (for $g$).
	Let 		$G : TX^*\to TX$ be defined by 
	\begin{eqnarray*}
\forall \alpha\in TX^*,\  \langle G(\alpha),\cdot \rangle_g = \alpha,
\end{eqnarray*}
where $g$ denotes the Riemannian metric. Since $\nabla s (x): T_xX\to E_x$ is onto, there exist bilinear forms $(k_i)_{i=1,\cdots, r}$ on $T_xZ$ such that for any local pair of tangent vector fields $V,W \in TZ$, the Weingarten operator $\sigma$ satisfies 
\begin{equation}\label{phit} \sigma(V,W)=(\nabla_V W)^\perp 
	= \sum_{i=1}^r k_i(V,W) G \langle \nabla s,t_i\rangle_h,
	\end{equation}
	so that 
	\begin{equation}\label{phir}
	\left(\langle(\nabla_V W)^\perp, G \langle\nabla s,t_i\rangle_h\rangle_g\right)_{i\in \{1, \cdots, r\}}
	= \left(\sum_{i=1}^r \Phi_{ij}k_i(V,W)\right)_{1\leq j\leq r} ,
	\end{equation}
	where 
	$$ \Phi = \left(\langle G \langle\nabla s,t_i\rangle_h G \langle \nabla s,t_j\rangle_h\rangle_g\right)_{1\leq i,j\leq r}$$
Using that the connection $\nabla$ is metric, (\ref{phir}) is equivalent to
	$$ (k_i(V,W))_{i=1,\cdots, r} = -\Phi^{-1}\left((\langle W, \nabla_V G \langle\nabla s,t_j\rangle_h\rangle_g)_{j=1, \cdots, r}\right).$$
	Since $\nabla G = G \nabla$, (\ref{phit}) becomes
	$$ \sigma(V,W)=-
	\sum_{i=1}^r (\Phi^{-1}\left(\langle W,  G \nabla_V \langle \nabla s,t_j\rangle_h\rangle_g)_{j\in \{1, \cdots, r\}}\rangle\right)_i G \langle \nabla s,t_i\rangle_h
	$$
and since $\nabla$ is metric, for any $j=1,\cdots, r$, 
\begin{eqnarray*}
	\langle W,  G \nabla_V \langle \nabla s,t_j\rangle_h\rangle_g &=& 
d_V \langle\nabla_W s,t_j\rangle_h- \langle\nabla_{\nabla_V W} s,t_j\rangle_h\\
	& =& 
	\langle\nabla_V\nabla_W s,t_j\rangle_h+	\langle\nabla_W s,\nabla_Vt_j\rangle_h
	- \langle \nabla_{\nabla_V W} s,t_j\rangle_h \\
 &=& \langle\nabla^2_{VW}s, t_j\rangle_h.
	\end{eqnarray*}
Here, we used that $\nabla_W s=0$ for any tangent vector $W\in T_xZ(s).$
Now,
		$$ \Phi = \left(\langle \langle\nabla s,t_i\rangle_h, \langle\nabla s,t_j\rangle_h\rangle_{g^{*}}\right)_{1\leq i,j\leq r}$$
		where $g^*$ denotes the scalar product on $TM^*$ associated to $g$.
Then, $\Phi$ is the matrix of the morphism $$\nabla s G \nabla s^* \in \mathcal L(E^*,E)$$
in the orthonormal basis $(t_j)_{1, \cdots, r}$ of $E$ and its dual basis $(t_j^*)_{j}$ of $E^*$ and $\nabla s^* : E^*\to TM^*$ is the adjoint of $\nabla s$. 
Hence,
	\begin{eqnarray*}
		\sigma(V,W) & =&
		G \langle \nabla s, ( \nabla s G \nabla s^*)^{-1}(\nabla^2_{V,W}s)\rangle_h
	\end{eqnarray*}
so that after a short computation we obtain:
\begin{equation}\label{uranus}
\langle \sigma(V,W), \sigma(X,Y)\rangle_g  =
( \nabla s G \nabla s^*)^{-1}(\nabla^2_{V,W}s) (\nabla^2_{X,Y} s).
	\end{equation}
\end{preuve}
\begin{corollary}\label{coroK}
	Under the hypotheses of Proposition~\ref{Marilyn}, for any $x\in Z(s)$,
	the sectional curvature $K$ of $(Z(s),g_{|Z(s)})$ at $x$ satisfies 
	$$ \max_{P\in \grass_{\R}(2,T_xX)}|K(x,P)-K^X(x,P)|\leq 
	2	\| (\nabla s G \nabla s^*)^{-1}\| \|\nabla^2 s\|^2.
	$$
	where $T_xM$ and $E_x$ are identified with their dual through their metric and the right-hand-side
	is evaluated at $x$.
\end{corollary}
\begin{preuve}
This is a straightforward consequence of Proposition~\ref{Marilyn} and the definition~(\ref{cas}) of the sectional curvature. 
\end{preuve}
\begin{lemma}\label{transfert} Let $1\leq r \leq n$ be two integers, $(F,g)$ and $(E,h)$ be two Hermitian spaces of finite dimensions $n$ and $r$, and $f: F \to E$ be of rank $r$. Then 
	$$ 
	\| (f  f^*)^{-1}\|\leq 
	\left(\min_{g\in \mathcal L_{sing} (F,E)} \| f-g\|^2\right)^{-1},
	$$
	where $E$ and $F$ are identified with their dual through their metric. 
\end{lemma}
\begin{preuve}[ of Lemma~\ref{transfert}] Let $J$ be the  orthogonal complement of $\ker f\subset F$. Then, 	$$ \min_{g\in \mathcal L_{sing} (F,E)} \| f-g\|^2 =
	\min_{g\in \mathcal L_{sing} (J,E)} \| f_{|J}-g\|^2.$$
	Let $A\in M_r(\C)$ be the matrix of $f_{|J}$ in two orthonormal basis of $F$ and $E$ respectively. 
	By the polar decomposition, 
	$	 A =  PU $, where $U$ is a unitary matrix and $P$ is a positive semi-definite Hermitian matrix. Clearly, 
	$$\min_{g\in \mathcal L_{sing} (J,E)} \| f_{|J}-g\|^2=
	\min_{B\in M_{r,sing} (\C)} \| P-B\|^2.$$
	Now, using an orthonormal basis of eigenvectors of $P$, we find that
	$$\min_{B\in M_{r,sing} (\C)} \| P-B\|^2= (\min \spec P)^2,$$
	where $\spec$ denotes the spectrum of the morphism.
	Now the matrix of $ff^*$ is $AA^* =P^2$ so that
	$$ \| (ff^*)^{-1} \|\leq (\min_{g\in \mathcal L_{sing} (F,E)} \| f-g\|^2)^{-1} .$$
\end{preuve}
Note that the proof shows that $$\min_{g\in \mathcal L_{sing} (F,E)} \| f-g\|^2
=\min \spec ff^*.$$
\begin{preuve}[ of Theorem~\ref{theoremcurvature}]
Let $(a_d)_d$ be a sequence converging to 0, $A$ be the event 
$$ A =\{s\in S^d, \ 
\dist_{\Delta_d} (s)
\geq a_d d^{-2n}\}$$
and for $k=\{0,1,2\}$, $B_k$ be the event
$$ B_k = \left\{s\in S^d, \ C
\sup_X \|\nabla^k s\|\leq  C\sqrt{d^k \log d}
\right\}$$
and $B =  A\cap B_0 \cap B_1 \cap B_2.$
By Proposition~\ref{burg} and Theorem~\ref{borne}, 
there exists $C$ depending on the sequence $(a_d)_d$ such that 
$$\forall d\geq 1, \  \mu_d (B) \geq 1- C(a_d+\frac{1}d).$$ 
From now on, assume that $s\in B$. 
By Proposition~\ref{raf2}, since  $s\in A$, 
there exists $C'>0$ a universal constant such that 
\begin{equation}
\min_{M\in \mathcal L_{sing}(T_xX, (E\otimes L^d)_x)}\|\nabla s(x)-M\|_{h^d}
\geq 
(a_d d^{-2n})d^{\frac{n+1}2}.
\end{equation}
By Proposition~\ref{curv} and Lemma~\ref{transfert}, for any $x\in X$,
	\begin{eqnarray*}\max_{P\in \grass_{\R}(2,T_xX)}|K(x,P)-K^X(x,P)|&\leq &
	C
	\| (\nabla s  \nabla s^*)^{-1}\| \|\nabla^2 s\|^2 \\
	& \leq & \frac{C}{a_d} d^{{\frac{3n-1}{2}}}d^2 \log d =\frac{C}{a_d} d^{\frac{3}2(n+1)} \log d,
\end{eqnarray*}
hence the conclusion.
	\end{preuve}

\subsection{Injectivity radius}

\begin{theorem}\label{ehrlich}(\cite[p. 156]{ehrlich})
	Let $(M,g)$ be a compact smooth Riemannian manifold. Assume that 	there exists $k>0$, such that 
	the sectional curvatures are bounded above by $k$. Then
	$$\inj (M,g)\geq \min\{\frac{\pi}{\sqrt k},	\frac12 \syst(M,g)\}.$$
\end{theorem}
\begin{preuve}[ of Theorem~\ref{theoremradius}]
This is a direct consequence of Theorem~\ref{sy}, Theorem~\ref{theoremcurvature} and Theorem~\ref{ehrlich}.
\end{preuve}

\subsection{Spectral gap}

\begin{theorem}\label{eig}(\cite[Theorem 1.2.1]{hassannezhad2016eigenvalue})
	 Let $(M, g)$ be a compact smooth Riemannian $n$-manifold. Assume that there exists $\kappa>0$ such that 
	 $$ \forall x\in M,\  \forall X\in \mathbb S T_xM, \ \ricci (X,X) \geq -(n-1)\kappa.$$ Then there exist $C>0$  depending
	on the dimension $n$ of $M$ only, such that
$$	\lambda_1\geq C^{1+\diam(M,g)\sqrt{\kappa} }\diam(M,g)^{-2}.$$
\end{theorem}
Recall that the Ricci curvature at a point $x$ in the direction $X\in T_xM$ is defined by
$$ \ricci (X,X)= \sum_{i=1}^n K(X,e_i),$$
where $(e_i)_{i=1,\cdots, n}$ is an orthonormal basis of $T_xM$. 
is the average of the sectional curvatures of all planes in $T_xM$. 

Hence, we need bounds for the diameter. In~\cite{feng1999diameter}, 
the authors proved that complex planar curves of degree $d$ have a diameter
less than $Cd^4$. 
We extend this theorem in our more general setting and with a better bound:
\begin{theorem}\label{dia}Under the hypotheses of Theorem~\ref{sy},
	There exists $C>0$ such that 
	$$ \forall d\gg 1, \forall s\in H^0(X, E\otimes L^d ),\ 
	\diam(Z(s))\leq Cd^{3r}.$$
\end{theorem}
The proof, which follows the lines of the one of~\cite{feng1999diameter}, is given in the next section.
\begin{preuve}[ of Theorem~\ref{theoremeigenvalue}]
	This is a direct consequence of Theorem~\ref{eig}, 
	Theorem~\ref{dia} and Theorem~\ref{theoremcurvature}.	
\end{preuve}
\section{Proof of the diameter's bound}\label{section: diameter}
In this section we prove Theorem~\ref{theorem: diameter}. Our proof is close to the one given by Feng and Schumacher~\cite{feng1999diameter}. We provide a better bound for plane algebraic curves, that is $d^3$ instead of $d^4$, and we adapt their proof to our general setting to prove the bound in any codimension and in any projective variety. 

Let $X$ be a complex projective manifold of dimension $n$ equipped with a K\"ahler metric $\omega$. Let us fix once for all an embedding $X\subset \C P^N$. Let $\omega_{\C P^N}$ be the Fubini-Study metric of $\C P^N$. 
Given a line $\C P^1\subset \C P^N$, the linear projection $\C P^N\dashrightarrow \C P^1$ defines a pencil of hyperplanes.
\begin{lemma}\label{lemma: pencils}
	There exists $\pi_1:\C P^N\dashrightarrow \C P^1,\dots, \pi_m:\C P^N\dashrightarrow \C P^1$ linear projections and $c_1>0$ such that $$\omega_{\C P^N}\leq c_1 \sum_{i=1}^m\pi_i^*\omega_{\C P^1}.$$
\end{lemma}
\begin{preuve} Let $x\in\C P^N$. By a direct computation, one can find finitely many linear projections $\pi_1:\C P^N\dashrightarrow \C P^1,\dots, \pi_k:\C P^N\dashrightarrow \C P^1$ such that $\sum_{i=1}^k\pi_i^*\omega_{\C P^1}$ is strictly positive at $x$. Being strictly positive is an open property, so $\sum_{i=1}^k\pi_i^*\omega_{\C P^1}$ is strictly positive on a neighborhood of $x$. Remark that, outside this neighborhood, $\sum_{i=1}^k\pi_i^*\omega_{\C P^1}$ is positive (meaning $\geq 0$) and that the sum of a strictly positive form with a positive one is strictly positive. 
	By compacity of $\C P^N$, one can than find $\pi_1:\C P^N\dashrightarrow \C P^1,\dots, \pi_m:\C P^N\dashrightarrow \C P^1$ such that $\sum_{i=1}^m\pi_i^*\omega_{\C P^1}$ is strictly positive. Then, again by compacity of $\C P^N$, there exists $c_1\in \R^*$ so that $\omega_{\C P^N}\leq c_1 \sum_{i=1}^m\pi_i^*\omega_{\C P^1}$.
\end{preuve}

From now on, we fix such pencils $\pi_1:\C P^N\dashrightarrow \C P^1,\dots, \pi_m:\C P^N\dashrightarrow \C P^1$ given by the previous lemma. 
Remark that if we restrict these pencils to $X$, we obtain pencils of hypersurfaces in $X$. By construction, all these hypersurfaces are hyperplane sections, all lying in the same cohomology  class $H\in H^2(X,\mathbb{Z})$.
Let $S\subset X$ be a smooth complex submanifold of dimension $r+1$, obtained by intersecting $X$ with a generic $\C P^{N-n+r+1}\subset \C P^N$. Then, for any $s\in H^0(X,E\otimes L^d)$, the intersection $S\cap Z(s)$ defines a complex curve $C(s)$ in $S$. We denote by $$u_i:C(s)\rightarrow \C P^1$$ the restriction of the map $\pi_i$ to $C(s)$.

\begin{lemma}\label{lemma: degree of branched} Using the previous notations, any linear projection $\pi:\C P^N\dashrightarrow \C P^1$ restricted to $C(s)$ defines a branched covering $u:C(s)\rightarrow \C P^1$ of degree $$ \sum_{i=0}^rH^{n-r}\cap c_{r-i}(E)\cap c_1(L)^id^i,$$ where $H$ is the class of the hyperplane section given by the embedding $X\subset \C P^N$ and $c_{\ell}(E)$ and  $c_1(L)$ are the Chern classes of $E$ and $L$.
\end{lemma}
\begin{preuve}
	Remark that the degree of a branched covering $u:C\rightarrow \C P^1$ equals the degree of the line bundle $u^*\mathcal{O}_{\C P^1}(1)$. In our case, the line bundle $u^*\mathcal{O}_{\C P^1}(1)$ equals the restriction to $C(s)$ of the line bundle $\mathcal{O}_{\C P^N}(1)$, so that his degree equals $ [C(s)]\cap c_1\left(\mathcal{O}_{\C P^N}(1)|_{C(s)}\right)$. Using the naturality of the cap product and that the class of $C(s)$ inside $X$ equals $[S]\cap c_r(E\otimes L^d)$ we see that the degree  $ [C(s)]\cap c_1\mathcal{O}_{\C P^N}(1)|_{C(s)}$ equals $$([S]\cap H\cap c_r(E\otimes L^d))=[S]\cap H\cap \sum_{i=0}^rc_{r-i}(E)\cap c_1(L^d)^i,$$ where $[S]\in H^{2(r+1)}(X,\mathbb{Z})$ is the class of $S$ in $X$ and where in the equality we used \cite[Example 3.2.2]{fulton}. The result then follows from the fact that $[S]=H^{n-r-1}$ and from $c_1(L^d)=dc_1(L)$.
\end{preuve}

\begin{lemma}\label{lemma: length of geodesic} Using the previous notations, let $\ell\subset \C P^1$ be a geodesic. Then $$\mathrm{length}\left(u_i^{-1}(\gamma)\right)\leq c_2d^{2r},$$ where $c_2$ is a constant depending only on the embedding $X\subset \C P^N$ and on $L$.
\end{lemma}
\begin{preuve} By Lemma \ref{lemma: degree of branched}, the degree of the branched covering $u_i: C(s)\rightarrow \C P^1$ is bounded by $c_3d^r$, where $c_3$ only depends on the embedding $X\subset \C P^N$ and on $L$. By the proof of ~\cite[Lemma 2]{feng1999diameter}, the curve $u_j\left(u_i^{-1}(\gamma)\right)$ is a  real algebraic curve of degree $c_4d^{2r}$, where $c_4$ only depends on the embedding $X\subset \C P^N$ and on $L$. 
	By Lemma 3 of Feng-Schumacher, if $\gamma$ be a real algebraic curve of degree $k$ inside $\C P^1$, then $\mathrm{length}_{\C P^1}(\gamma)\leq 2\pi k$.  Moreover, by Lemma \ref{lemma: pencils}, we have $\mathrm{length}_{C(s)}(u_i^{-1}(\gamma))\leq c_1\sum_{j=1}^m\mathrm{length}_{\C P^1}\left(u_j(u_i^{-1}(\gamma))\right)$, where $c_1$ does not depend on $C(s)$.  This proves the lemma.
\end{preuve}
\begin{preuve}[ of Theorem \ref{theorem: diameter}]
	Given $s\in H^0(X,E\otimes L^d)$, we want to bound from above the diameter of $(Z(s),\omega|_{Z(s)})$. Let us fix once for all an embedding $X\subset \C P^N$ and  $\pi_1:\C P^N\dashrightarrow \C P^1,\dots, \pi_m:\C P^N\dashrightarrow \C P^1$ as in Lemma \ref{lemma: pencils}.
	As $\omega\leq c_1\omega_{\C P^N}|_{X}$ for some $c_1>0$, it is enough to bound from above the diameter of $(Z(s),\omega_{\C P^N}|_{Z(s)})$.
	In order to bound from above the diameter of $Z(s)$, we will bound from above the distance between any pair of points $p,q\in Z(s)$. To do this, we will consider a generic $\C P^{N-n+r+1}$ passing through $p$ et $q$, which defines a $(r+1)$--dimensional complex submanifold  $S=X\cap \C P^{N-n+r+1}$ together with a curve $C(s)=Z(s)\cap \C P^{N-n+r+1}$. We have the inequality $\mathrm{dist}_{Z(s)}(p,q)\leq \mathrm{dist}_{C(s)}(p,q)$, so that the result will follow from the inequality 
	\begin{equation}\label{inequality for the distance}
	\mathrm{dist}_{C(s)}(p,q)\leq cd^{3r}.
	\end{equation}
	In order to prove \eqref{inequality for the distance}, we will explicitly construct a path between $p$ and $q$  whose length is bounded from above by $cd^3$. Here is the construction of such path. Let 
	$$u:=u_1:C(s)\rightarrow \C P^1$$ be the restriction of $\pi_1:\C P^N\dashrightarrow \C P^1$ to $C(s)$.  By Lemma \ref{lemma: degree of branched}, the degree of $u$ is bounded from above by $D=c_3d^r$, where $c_3$ only depends on the embedding $X\subset \C P^N$ and on $L$ and $E$. We can assume that $u$ is a simple branched covering, that is all the branched points are simple (i.e. of multiplicity $2$) that $p$ and $q$ are not branched points. Let denote by $p_j\in C(s)$, $j\in{1,\dots,\ell}$, the branched points. We denote by $q_j=u(p_j)\in \C P^1$ the critical values. Eventually after an arbitrarily small perturbation of $u$, we can suppose that  the $q_j$'s are all distinct and that no three of them are contained in a closed geodesic.

	We choose an auxiliary point 
	$$x\in \displaystyle\C P^1\setminus \{q_1, \cdots, q_\ell\}$$

and for any $j\in \{1,\cdots, \ell\}$  we consider an arc of geodesic $S_j$ from $x$ to $q_j$. Then, $\displaystyle\C P^1\setminus \cup_{j}S_j$ is contractible and then the cover $$u:u^{-1}(\displaystyle\C P^1\setminus \cup_{j}S_j)\rightarrow \displaystyle\C P^1\setminus \cup_{j}S_j$$ is the degree $D$ trivial cover, consisting of $D$ copies $F_1,\dots,F_D$ with $F_i\simeq \C P^1\setminus \cup_{j}S_j$, called the sheets of the covering. We denote by 
$S^{(i)}_j$ and by $x_i$ the copies of 
$S_j$ and of $x$ appearing in the $i$--th sheet. 
	
	As the branched covering $u:C(s)\rightarrow \C P^1$ has only simple ramifications, any branched point $p_j$ is contained in the closure of precisely two sheets. Any two sheets whose closures share a branched point are called \emph{adjacent}. If $a$ and $b$ are two points lying in two adjacent sheets, then there is a path in $C(s)$ of length bounded by $c_4d^2$ that goes from $a$ to $b$. Such path is constructed as follows. Let $F_{i(a)}$ and $F_{i(b)}$ be the sheets containing $a$ and $b$ and let $p_{j(a,b)}$ be the branched point connecting $\bar{F}_{i(a)}$ and $\bar{F}_{i(b)}$. We can that go from $a$ to $x_{i(a)}$ following a path which is contained in $u^{-1}(\gamma_1)$, where $\gamma_1$ is a geodesic going from $u(a)$ to $x$. Then we can go from $x_{i(a)}$ to $p_{ab}$ following $S^{(i(a))}_{j(a,b)}$, which is, by construction, contained in  $u^{-1}(S_{j(a,b)})$. As the point $p_{j(a,b)}$ lies in $\bar{F}_{i(b)}$, we can follow backwards $S^{(i(b))}_{j(a,b)}$ and go from  $p_{j(a,b)}$ to $x_{i(b)}$. Such path is again contained in $u^{-1}(S_{j(a,b)})$. Finally, we can go from $x_b$ to $b$  following a path which is contained in $u^{-1}(\gamma_2)$, where $\gamma_2$ is a geodesic going from $u(b)$ to $x$. By Lemma \ref{lemma: length of geodesic}, each of these four paths has length bounded from above by $c_2d^{2r}$, where $c_2$ is a constant depending only on the embedding $X\subset \C P^N$ and on $L$, so the distance between $a$ and $b$ is bounded from above by $4c_2d^{2r}$.
	
	We have then proved that any pair of points lying in two adjacent sheets can be joint by a path of length smaller than $c_4d^{2r}$. This implies that any pair of points lying in two $k$--adjacent sheets (two sheets are called $k$--adjacent if one can pass from one to the other by passing through at most $k$ branched points) can be joint by a curve of length smaller than $c_4d^{2r}k$. As every pair of sheets is $(D-1)$--adjacent (recall that $D$ is the degree of the covering) and as $D=c_3d^{r}$ (see Lemma \ref{lemma: degree of branched}), this implies \eqref{inequality for the distance}. Hence the theorem. 
\end{preuve}

\bibliographystyle{amsplain}
\bibliography{metricaspects.bib}

\noindent
Michele Ancona \\
Laboratoire J.A. Dieudonn\'e\\
UMR CNRS 7351\\
Universit\'e C\^ote d'Azur, Parc Valrose\\
06108 Nice, Cedex 2, France\\

\noindent
Damien Gayet\\
 Univ. Grenoble Alpes, Institut Fourier \\
F-38000 Grenoble, France \\
CNRS UMR 5208  \\
CNRS, IF, F-38000 Grenoble, France
\end{document}